\documentclass[preprintnumbers,amsmath,amssymb]{revtex4}
\usepackage{graphicx}
\usepackage{graphicx}
\usepackage{bm}
\def\dir{}
\def\figs{}

\begin{document}

\title{Truels, or the survival of the weakest}

\author{Pau Amengual}
 \homepage{http://www.imedea.uib.es/~pau/}
 \email{pau@imedea.uib.es}
\author{Ra\'ul Toral}%
 \email{raul@imedea.uib.es}
\affiliation{%
Instituto Mediterr\'aneo de Estudios Avanzados (IMEDEA)\\\
CSIC-University of the Balearic Islands\\
Ed. Mateu Orfila, Campus UIB, E-07122 Palma de Mallorca, Spain
}%

\maketitle

\section{Introduction}

We are all familiar with the concept of a ``duel", a fight between two people whose main goal is to eliminate the opponent. The goal can be very aggressive, as in the duels involving offended lovers in the romantic times and the duels between the sheriff and the gunman in the far-west, or it could end in both players enjoying a beer, such as in a darts game between two friends (with the loser paying the beer). Using a non-bloodily point of view, a duel is simply a game with two players, each one having a probability of winning the game. The players have an intrinsic ``marksmanship" or ability associated to their performance in the game.  Of course, one could favor the worst player, the one with the worst marksmanship, with some advantage, as some extra bullets or some extra points, but in a normal, not biased situation, the common sense tells us that the higher the marksmanship, the higher the probability of winning the duel. The mathematical treatment of this game confirms our simple expectations. If there are more than two players, a series of duels might be needed in order to determine the absolute winners. This is the case, for instance, of a basketball league\footnote{A football league is much more interesting, of course, but a football game can end up in a draw, a situation not considered in the duels.}, where the teams are paired in all possible ways. The winner of the league is the team that wins more games. 

A ``truel" is a generalization of a duel involving  three players and with some rules that we now spell out in detail. First, it is assumed that the marksmanship of a player does not depend on the opponent. For instance, the probability that a person kills the opponent in the romantic duel does not depend on who is the opponent. Second,  ``truelists" shoot sequentially, one after another. If a truelist misses the shot, then it is the turn of another truelist. Third, each truelist has the right to chose the person to aim at. More specifically, the steps are:
 
\begin{enumerate}
\item One of the remaining truelists is chosen at random.
\item The chosen truelist selects an opponent and, with a certain probability --the marksmanship-- eliminates that opponent from the game.
\item Whatever the result, steps 1 and 2 are repeated until there is only one survivor, the winner of the game.
\end{enumerate}

These rules define what is known as the {\sl random truel}. A paradoxical result appears, as the player with the highest marksmanship does not necessarily possess the highest survival (or winning) probability. This paradoxical result was already mentioned in the early literature on truels~\cite{kb97.1,s82.1}. The truel is an example of a situation in which ``the fittest does not necessarily survive"\cite{s54.1}. We all know cases in which the more appropriate person is not the one who gets the job or the one who reaches to the top of an organization, and the truels present a framework in which this counterintuitive outcome can be rigorously derived. It seems that these games were introduced for the first time in~\cite{k46.1}, although the name \textit{truel} was coined later by Shubik in the 1960's. 

Besides the potential applications, truels are also interesting from the pedagogical point of view since they illustrate in a particularly clear way some of the concepts of game theory. In particular, the Nash equilibrium can be understood as the best action players can take to optimize their probability of winning the game. The ``action" here corresponds to the strategy when choosing the opponent. Intuition tells us that the best one can do given the above rules  is to shot at the best of the remaining players, and this result is confirmed by the mathematical analysis of the game. This ``best opponent strategy" explains the paradoxical result, since it implies that the best player is the target of the other two, decreasing accordingly that player's probability of winning. 

This intuitive explanation of the paradoxical result needs to be confirmed by a rigorous mathematical treatment. This has been done previously and extensively using the natural language of game theory\cite{k72.1,k75.1,k77.1}. We have revisited\cite{ta05.1} this analysis from the point of view of stochastic processes using discrete--time Markov chains~\cite{k73.1} with three absorbing states. The states of the chain are determined by the players who are still in the game, and the absorbing states correspond to those states where there is only one player remaining. By means of this technique we can calculate the probability of the system ending in one of those three absorbing states, that is, the probability of winning for every player.
This point of view, besides reproducing the main paradoxical result, provides a framework for the analysis of the games with which some communities (such as the physicists) feel more comfortable. In fact, Physics has contributed recently to the field of paradoxical games with the so called Parrondo's games\cite{ha99.1,ha02.1,tam03.1}. These games exemplify the situation in which a combination of losing games results in a winning game. They are inspired by the ratchet and pawl which is a pedagogical device invented by Smoluchovsky\cite{ms12.1}, and largely studied by Feymann\cite{fls63.1,pe96.1}, in order to illustrate the second law of thermodynamics, the one that poses a limit to the amount of energy we can extract from a given source.

We  find in the literature other models similar to the truel game that present also counterintuitive results, like for instance the \textit{rock--scissors--paper} game. This game has been applied to population dynamics~\cite{sl96.1,fa01.1} consisting on a system with three species that interact with each other creating a competitive loop (recall that in the  \textit{rock--scissors--paper} game a rock beats a pair of scissors, scissors beat a sheet of paper and paper beats a rock). The paradoxical effect in this model is that the least competitive species might be the one with the largest population and, when there are oscillations in a finite population, to be the least likely to die out. This game has also been applied to a voter model\cite{t93.1,t95.1} obtaining again a paradoxical result, namely, an initial damage and suppression of one candidate may later lead to an enhancement of the same candidate.

Several modifications of the basic rules of truels have been proposed. For instance, one could bias the game by allowing the worst player to shoot in the first place, with the best a priori player shooting last. This modification gives the worst player an initial advantage. It is of no surprise that this ``sequential rule" increases the probability that the worst player wins the game. It is surprising, though, that in order to increase the survival probability, the worst player might need to fail the shots {\sl on purpose}. A yet different version, not considered in this paper, is that of the {\sl simultaneous truel} where the three players shoot at the same time. 

Other rules could also be modified. For instance, one can consider that the number of rounds is finite or infinite, or the ammunition to be limited or unlimited, and these modifications do modify the probabilities of the different outcomes~\cite{k72.1,k75.1,k77.1}. Another modification leads to the \textit{cooperative truels}~\cite{bbk02.1}. They are characterized by the appearance of cooperations in which different players set a common target and improve, by means of such a coalition, their survival probability. 

In this paper we review some of the main results  obtained in this field of truels.  We limit ourselves to the random and sequential truels in which players use their best possible strategy with no coalitions. First, in order to introduce the notation and the main concepts, we will present the simpler case of duels. We then consider how the appearance of a third player modifies substantially the analysis and players must decide on a certain strategy. The ``best strategy" in a sense to be precisely defined leads to the Nash equilibrium point as the best action that rational players can take. The use of this best strategy may result in the worst player winning the game. 

We have modified the random truel explained above and converted it into an opinion model. In this version each of the three players holds a different opinion on a given topic. Otherwise, the mechanics of the game is similar to the random truel. At each round one person is chosen randomly amongst the players and then that person tries to change the opinion of another player.  The game ends when all players share the same opinion.

We address next the question of who wins a ``truel league". We will see that, despite the paradoxical result mentioned above, still the distribution of winners is peaked around the players with the higher marksmanship for the random and opinion versions. In the sequential truel, however, the paradoxical result remains partially since the distribution of winners is peaked around the intermediate players.

If the rules of truels are extended from three to $N$ players, the paradoxical results shows up even more clearly since as $N$ increases it is more difficult for the player with the highest marksmanship to win the game. Finally, we consider the dynamics of the games in the case that players have a well defined spatial distribution in a given network of interactions.

\section{The duels}

In this game there are two players: a  \textit{good} player (A) and a \textit{bad} player (B). Their respective marksmanships will be denoted by \textit{a} and \textit{b} respectively, with $a>b$. Here, the strategy is obvious: shoot at the only opponent. It makes no sense to lose an opportunity to eliminate the opponent by shooting into the air.

In the random case, in which the next player to shoot is chosen randomly at each time step, the analytical study (using, for example, the Markov chain's tools) shows that the winning probability of a player is proportional to the marksmanship of that player. If $P_A$ (resp. $P_B$) are the probabilities of player A (resp. B) winning the game, we have $P_A=\frac{a}{a+b}$ and $P_B=\frac{b}{a+b}$. Therefore, in the random duel the best player has the largest survival probability. The same result holds if the first player to shoot is chosen randomly and then they shoot sequentially one after the other until only one remains.

One could overcome this result by allowing the worst player to shoot first and alternate turns afterwards. In this sequential duel, the analysis of the resulting winning probabilities shows that only when $b>\frac{a}{1+a}$ the survival probability of player B is greater than that of A. Therefore, in the sequential duel, player B can overcome the unfavorable situation of having a lower marksmanship by being the one shooting in first place.

\section{Strategies in truels. Equilibrium points}

If a third player enters the game, the previous situation of a duel is no longer simple. Now every player in the truel must consider all possible actions that other opponents may take and their corresponding outcomes. In this case, we must consider strategies and make use of some concepts of game theory. For concreteness, and without loss of generality, we consider that the third player C has the lowest marksmanship, $c$, such that $a>b>c$.

All players in the truel share the same goal: to be the only one survivor. This can be explicitly imposed through the inclusion of a ``payoff", a concept introduced in game theory and that corresponds to some sort of reward the player receives for achieving the goal. In order to maximize their payoff, players have to chose strategies that maximize their survival probability. When the three players are still in the game, a player has three possible strategies\footnote{We only consider ``pure strategies". It is also possible to select with a given probability one of the options. These ``mixed strategies" are not relevant in the case of truels.}: two correspond to choosing one of the two opponents and the third strategy is to shoot in the air (or missing the shot on purpose). If one of the three players has been removed from the game, we are in a duel situation and, as discussed before, the only strategy is to aim at the remaining opponent. We also assume that strategies adopted by the players are non--cooperative, in the sense that alliances or pacts between them are not allowed. The corresponding expressions for the survival probabilities can be found in references \cite{k75.1,ta05.1}, but they are too cumbersome to be reproduced here. 

We consider first  the random truel in which the player to shoot is chosen randomly amongst the remaining players. If players do not use any strategy and simply shoot randomly to any of the two opponents it is easy to see that the winning probabilities are proportional to the marksmanship. However, and following a basic premise of game theory~\cite{or97.1}, we assume that players are {\sl rational agents}, in the sense that they adopt the strategy that is a best response to the strategies chosen by the other players. This is nothing but the concept of \textit{Nash equilibrium}, which we now clarify with an example\footnote{The well known movie {\sl A beautiful mind}, besides telling us some details of the life of J.F. Nash, contains some other examples of practical applications of the concept of Nash equilibrium.}. In Table~\ref{table_probs} we present the different survival --or winning-- probabilities $P_A$, $P_B$ and $P_C$ of players A, B and C respectively for the different strategies adopted by the players when they play the random truel. These values are calculated considering that player A has a marksmanship $a=1$ (or $100\%$ of effectiveness), player B has $b=0.8$ ($80\%$ of effectiveness) and player C $c=0.5$ ($50\%$ of effectiveness). For simplicity, we do not include in this example the strategy of shooting into the air, an option which is not relevant for the random truel.

Let's start by looking in Table~\ref{table_probs} at the set of strategies CCB. The notation means that the first player A aims at  C, the second player B aims at C and the third player C aims at  B. In this case the player with the highest survival probability is A with $58\%$ percentage of winning, followed by player B with $34.8\%$ percentage and finally player C with a very low percentage of $7.2\%$. When analyzing this situation, player C concludes that it is better to change the strategy and, instead of aiming at B, set player A as the new target. This changes to row CCA of the table where the  survival probability of player C has increased up to $8.5\%$.

Player B now looks at this new situation and decides to change the strategy and set A as the new target. We then move to the CAA row, where $P_B$ increases from $48.1\%$ to $54.1\%$. Now it is the turn of player A who decides to change strategy and set B as a new target thus leading the the set BAA where $P_A$ has indeed increased from $24.2\%$ to $29.0\%$ .

At this point, players can no longer change unilaterally and increase their survival probability. The important point is that if we repeat the argument starting from any strategy and using any order in the players reasoning, we will reach --sooner or later-- the same strategy set: BAA.

The strategy BAA is known in the literature~\cite{kb97.1,k75.1} as the \textit{strongest--opponent strategy}, as all players aim at the opponent with the highest marksmanship: player A aims at B, and players B and C aim at A. This is the unique \textit{Nash equilibrium point} of the random truel, meaning that no player improves the survival probability by changing strategy, as long as the rest of players keep theirs. Therefore, this set corresponds to a local maximum of all survival probabilities of the player and it is the one that rational players will necessarily use. It is remarkable that for this strategy set BAA (and the particular set of marksmanships chosen for $a$, $b$ and $c$), the survival probability goes in reverse with the marksmanship. Player C is the one with the highest survival probability and player A has the lowest. This is a counterintuitive result, in the sense that one would naively expect player A to have the highest survival probability on the basis that player A has the highest marksmanship.

\begin{table}[htb!]
\centerline{\begin{tabular}{|c|c|c|c|}
\hline
Strategy & $\mathbf{P_A}$ & $\mathbf{P_B}$& $\mathbf{P_C}$ \\
\hline
{CCB} & 0.580 & 0.348 & 0.072 \\
{CCA} & 0.434 & 0.481 & 0.085 \\
{CAB} & 0.386 & 0.407 & 0.207 \\
{CAA} & 0.242 & 0.541 & 0.218 \\
{BCB} & 0.628 & 0.155 & 0.217 \\
{BCA} & 0.483 & 0.288 & 0.229 \\
{BAB} & 0.435 & 0.214 & 0.351 \\
{BAA} & 0.290 & 0.348 & 0.362 \\
\hline
\end{tabular}}
\caption{Survival probabilities $P_A$, $P_B$ and $P_C$ of players A, B and C respectively, for the different set of strategies adopted in the case of the random truel. The marksmanships for player A, B and C are respectively $a=1$, $b=0.8$ and $c=0.5$. The Nash equilibrium point corresponds to the last row, BAA.\label{table_probs}}
\end{table}

In the sequential truel,  the worst player has the advantage of being the first one to shoot. The sequence of shots is C-B-A and it repeats as necessary until there is only one player left. This case is more complex than the random truel, and a detailed analysis reveals the existence of two different equilibrium points depending on the actual values of the marksmanships $a$, $b$ and $c$. One of the equilibrium points corresponds to the \textit{strongest--opponent strategy} (as in the case of the random truel), whereas the other corresponds to the set BA$\emptyset$, meaning that player A sets as a target player B, player B sets A as a target and player C decides to shoot into the air. This latter set of strategies implies that it might better for player C to fail the shot and let players B and A kill each other. Player C would use the next  turn to try to eliminate the remaining player, becoming the winner of the truel.

\begin{figure}[hbt]
\begin{center}
\includegraphics[scale=0.4,draft=false]{\figs 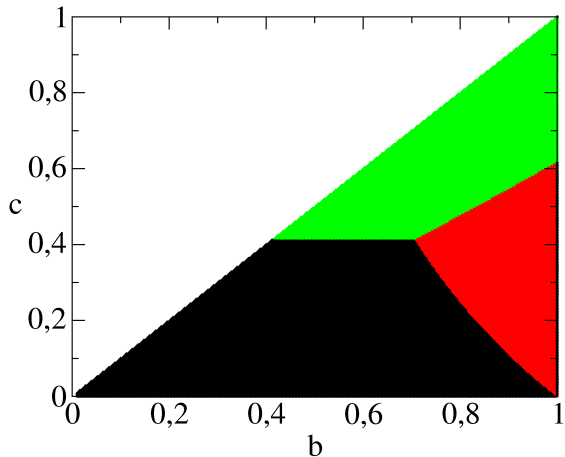}
\end{center}
\caption{The different colors (black for player A, red for player B and green for player C) indicate the regions in the parameter space $(b,c)$, setting $a=1$, in which a player has the highest probability of survival, i.e. is the favorite for winning the random truel game.\label{random}}
\end{figure}

Which player has the highest survival probability (i.e. who is the favorite to win the game) depends on the values of the marksmanships $a$, $b$ and $c$. In Fig.~\ref{random} we indicate the regions in the parameter space $(b,c)$, after setting $a=1$, in which each player of the random truel has the highest survival probability. Black corresponds to player A, green to B and red to C. This figure shows that it is possible for every player to be the favorite in the game, although the largest region corresponds to the best player A.

\begin{figure}[hbt]
\begin{center}
\includegraphics[scale=0.4,draft=false]{\figs 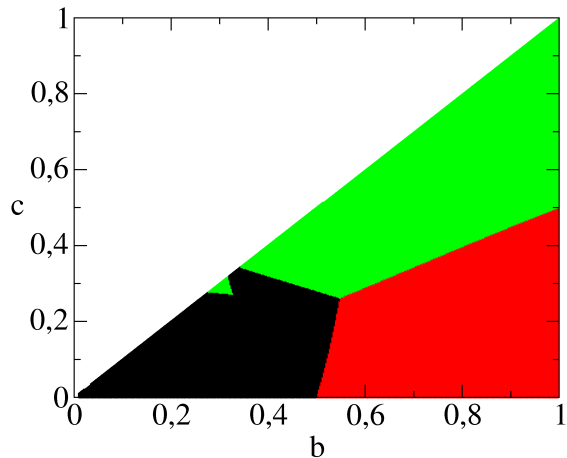}
\end{center}
\caption{Same as Fig. \ref{random} but for the case of the sequential truel.\label{sequential}}
\end{figure}

For the sequential truel, the situation we find is somewhat different. Due to the imposed firing order (C-B-A), player A is the last one to shoot. Therefore, the a priori advantageous situation given by a high marksmanship is partially lost. This is reflected in Fig.~\ref{sequential}, since  the region where player A is the favorite has decreased considerably compared to that of Fig.~\ref{random}. In fact, the a priori worst player C is the favorite in a larger number of occasions.  We explained previously that there were two equilibrium points in the sequential truel, BAA and BA$\emptyset$. The last one is the relevant in the small green region located in the black region seen in Fig.~\ref{sequential}.

\section{Truels as a model of opinion spreading}

We can give an interesting twist to a truel (as well as eliminating some of its aggressiveness) by interpreting it as a model of opinion spreading. In this case A, B and C are three opinions that  people can hold on a topic. In this version, players aim to convince, rather than eliminate, each other. Marksmanships $a$, $b$ and $c$ are now interpreted as the probability of convincing another person. The game ends when one of the opinions is shared by all the players.

\begin{figure}[hbt]
\begin{center}
\includegraphics[scale=0.4,draft=false]{\figs 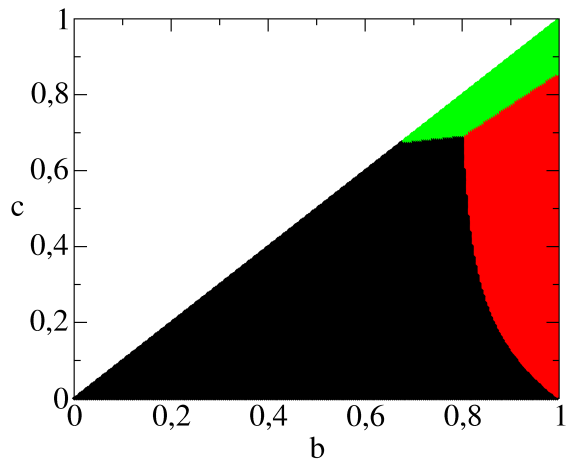}
\end{center}
\caption{Same as Fig. \ref{random} but for the case of the opinion version of the random truel.\label{opinion}}
\end{figure}

In the random version, the theoretical analysis shows that, independently on the values of the marksmanships, the only equilibrium point is the \textit{strongest--opponent strategy}. The same paradoxical result still applies as the opinion with the higher marksmanship does not necessarily need to be the one that survives. However, as shown in Fig.~\ref{opinion}, opinion A is the favorite to become the  majority opinion  for a larger region of  values of $b$ and $c$, again setting $a=1$. It is only for a relatively small region that opinion C is the favorite one. This overwhelming superiority of opinion A can be understood if we recall that in this model the total number of players remains constant throughout the game. Only the opinions held by the players change. So, once opinion A convinces either a player with opinion B or a player with opinion C, it is very likely that it will eventually become the majority opinion due to its high convincing probability. For instance, starting from the initial configuration where the three players held different opinions A, B and C, and for the set same values of $a$, $b$ and $c$ as before ($a = 1$, $b=0.8$ and $c=0.5$), there is a $38.6\%$ probability that opinion A becomes the majority opinion, and a $37.8\%$ and $23.4\%$ probability instead that opinions B and C become majority, respectively. 

\section{Distribution of winners}

Imagine now that we set a truel competition amongst a population of players.  The marksmanships of the population are uniformly distributed in the interval $(0,1)$. We set a league scheme. We form all possible triplets of players and have them play the random truel (we decide on a non lethal version in which losers are still able to play another game). The winner of the competition will be the player with the highest number of truels won. For whom would you bet? For the good players with high marksmanship or for the bad players with low marksmanship? In this case, one has to consider the average probability of winning for a player when playing against all sort of players. The mathematical analysis has been summarized in Fig. {\ref{f_various}. For the random truel and the opinion model, the histogram of winners has a maximum at a marksmanship of 100\%. In some sense, justice is restored in this truel competitions, since the best players are the ones who win the competition in more occasions. However, with the rules of the sequential truel in which the worst player is favored, it turns out  that the distribution is now peaked around an intermediate marksmanship of about 55\%. It's the triumph of the mediocrity!

\begin{figure}[hbt!]
\begin{center}
\includegraphics[scale=0.3,angle=-90,draft=false]{\figs 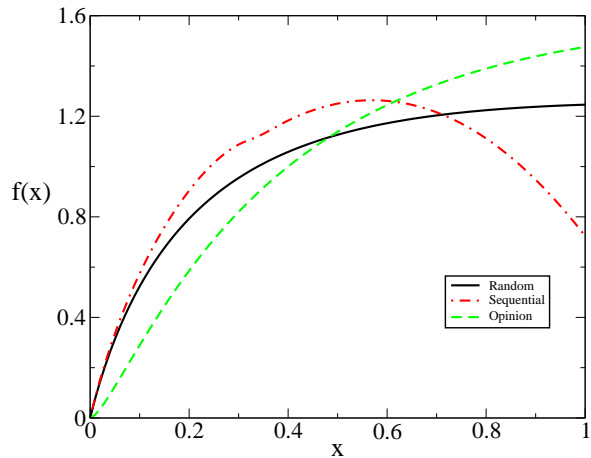}
\end{center}
\caption{Distribution of winners in a truel competition. The black (continuous) line corresponds to the random truel, the red (dot-dashed) to the sequential truel and the green (dashed) to the convincing opinion version of the model. Notice that in the case of sequential truel, the maximum occurs around $x=0.55$ while for the other two versions, the maximum is at $x=1$.\label{f_various}}.
\end{figure}

\section{Generalization to $N$ players : $N$--uels}

If we extend the rules of the truels to more than three players, the situation becomes very difficult to study analytically. However, it is rather easy to implement those rules in a computer program and visualize the results. Let us consider $N>3$ players in a random  $N$-uel. We use the same rules than the random truel as spelled out in the introduction chapter. When players use the strongest--opponent strategy, it is clear that $N-1$ guns will point at the best player. Hence the probability of survival of the best player is low and, moreover, decrease with increasing $N$. In the simulations we have considered sets of $N$ players whose marksmanship has been drawn from a uniform random distribution in the interval $(0,1)$. 

In Fig.~\ref{histo_N4} we show the histogram corresponding to the classification obtained when the game is played by $N=4$ players. The fourth classified corresponds to the distribution of players eliminated from the game in the first place, the third classified to  the ones eliminated in second place and so on. The distribution of the fourth classified shows that individuals eliminated firstly in the game are those with higher marksmanships. Indeed, the maximum is located at a marksmanship $x=1$, indicating that the better a player is the higher the probability of being eliminated first. Another aspect we can extract from this figure has to deal with the distribution of the first and second classifieds: these curves correspond to the case where there are only two players left in the game, i.e., to a duel. Therefore, it is more likely in this situation that players with lower marksmanships are eliminated firstly rather than those with higher marksmanships (that is the reason why the curve for the second classified presents a maximum in the origin). It is also worth mentioning that already for $4$ players the histogram associated to the first classified -- i.e., the winner of the $4$--uel -- presents a maximum at a value of $x<1$. This result implies that the best performing player does not correspond anymore to the player with the highest marksmanship, as it happened when $N=3$. Indeed, the optimum value is located around $0.45$. 

In Fig. ~\ref{various_N} we present the histogram of winners in a $N$-uel as a function of the number of players $N$. Accordingly with the previous discussion, it is evident that  the distribution is indeed progressively enhanced and shifted towards zero when $N$ is increased. In this limiting case, it is clear that we can talk about the ``survival of the weak ones".

\begin{figure}[hbt!]
\begin{center}
\includegraphics[scale=0.3,angle=-90,draft=false]{\figs 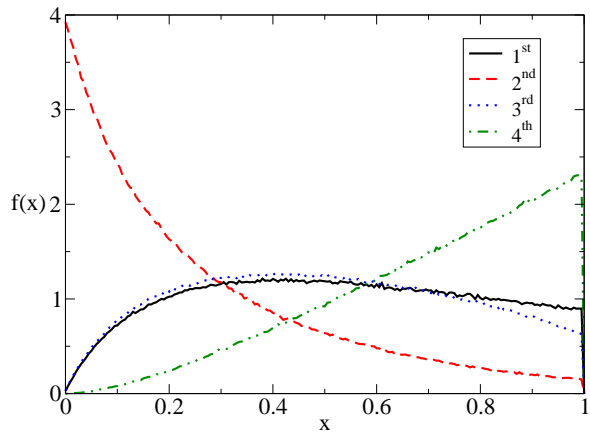}
\end{center}
\caption{Histograms of the first, second, third and classified corresponding to a random $N$-uel tournament for $N=4$ players.  \label{histo_N4}}
\end{figure}

\begin{figure}[hbt!]
\begin{center}
\includegraphics[scale=0.3,angle=-90,draft=false]{\figs 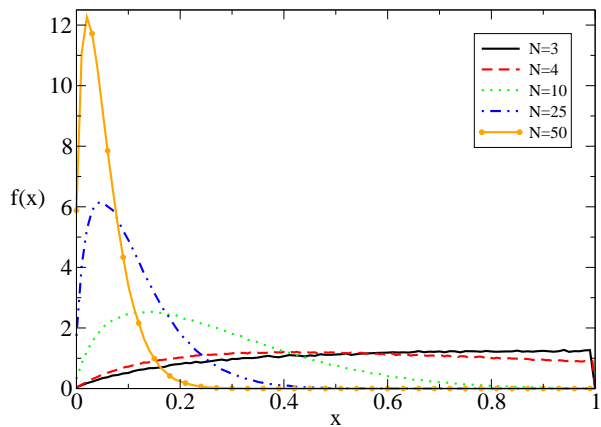}
\end{center}
\caption{Histograms of the winners of a random  $N$--uel tournament for different values of $N=3$, $4$, $10$, $25$ and $50$.\label{various_N}}
\end{figure}

\section{Truels with spatial dependence}

A natural step in the truels is the introduction of a spatial structure in the set of players. This reflects the fact that players do not interact with any other player, but only with those which are closer in some sense. Although one could devise some sort of social network of interactions\cite{ws98.1,dm02.1}, we consider here a simple  two--dimensional lattice with $N$ sites, each with  4 nearest neighbor links. The lattice is initialized by putting randomly on each site one player of groups A, B or C in the respective proportions $x_A$, $x_B$ and $x_C$, ($x_A+x_B+x_C=1$) and respective marksmanships $a$, $b$ and $c$. An important ingredient of this generalization is that players never shoot to a person of the same group. The rules of the random {\sl collective truel} are as follows:

\begin{enumerate}
\item One of the remaining players is chosen at random.
\item The chosen player selects randomly two players amongst the occupied neighbors sites and the three of them play a random truel. The losers of the truel are eliminated from the system. If the chosen player has only one neighbor, the two of them will play a duel with the loser being removed from the system. If no neighbors are left, the player will walk to a randomly chosen neighbor site. 
\item Steps 1 and 2 are repeated until all the survivors belong to the same group.
\end{enumerate}

In step 2, it is possible that some of the chosen  players belong to the same group.  In this case, they observe strictly the rule of no shooting between members of the same group. Accordingly, it could happen that there is more than one survivor of that game. In any event, players use the strongest--opponent strategy. If, for example, the three players in a truel belong to groups A, A and B, the two A players will aim at B, while B will aim to one of the two A (again chosen at random). The outcome of that particular situation could be either player B eliminating both A players or player B being eliminated by the two players A. Since the analytical treatment appears rather difficult, we present now the results coming from a direct numerical simulation of the aforementioned rules. We use throughout this section the values $a=1$, $b=0.8$, $c=0.5$ for the marksmanships.

In Fig.~\ref{evolution_random} we show some snapshots concerning different stages of a simulation carried out for the random truel. The initial population proportion was $x_A=0.3$, $x_B=0.3$ and $x_C=0.4$. Note that in the early stages of the simulation the populations of groups B and C diminish, whereas group A resists and eventually becomes the winner of the collective truel.

\begin{figure}[hbt!]
\begin{center}
\includegraphics[scale=0.22,angle=-90,draft=false]{\figs 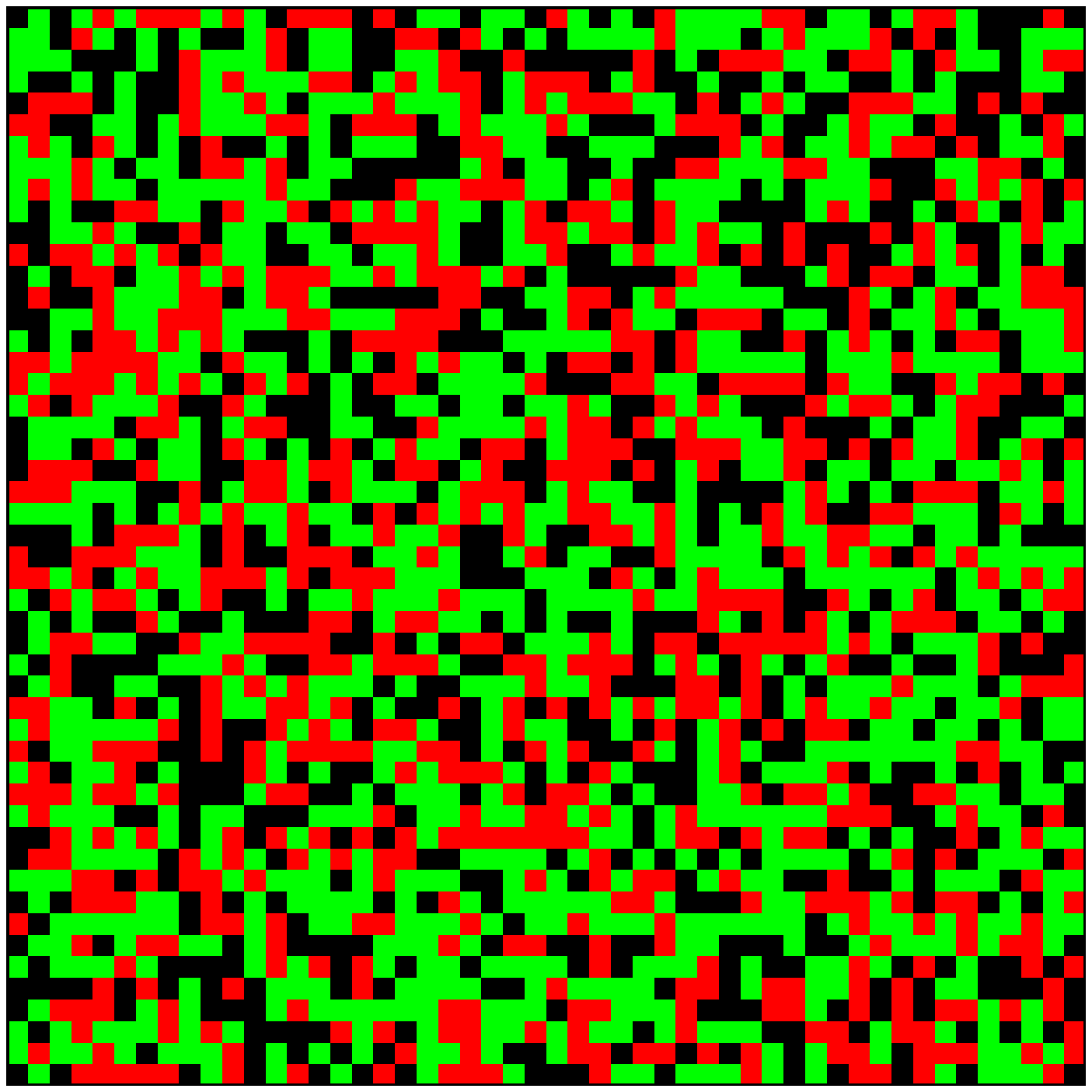}
\includegraphics[scale=0.22,angle=-90,draft=false]{\figs 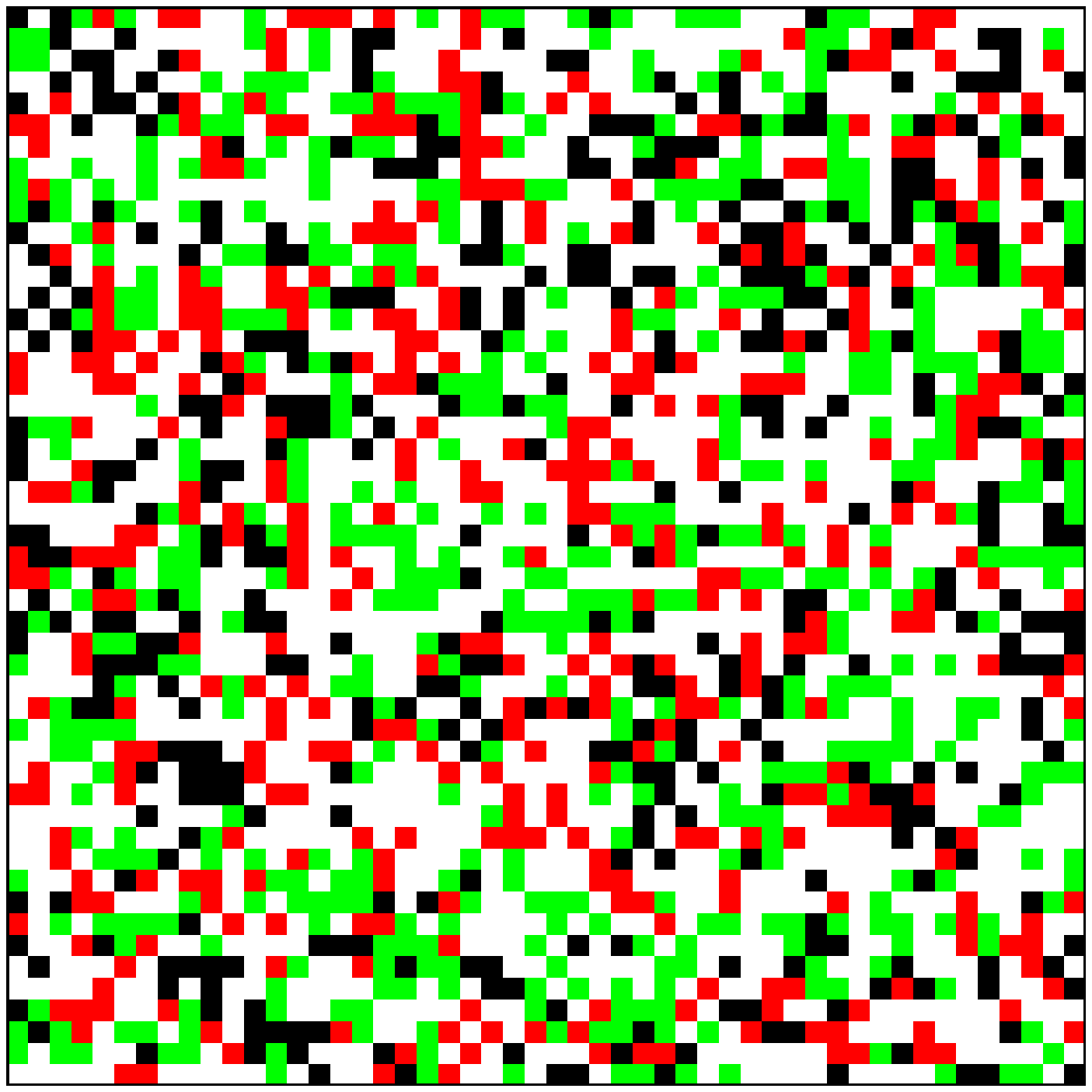}
\includegraphics[scale=0.22,angle=-90,draft=false]{\figs 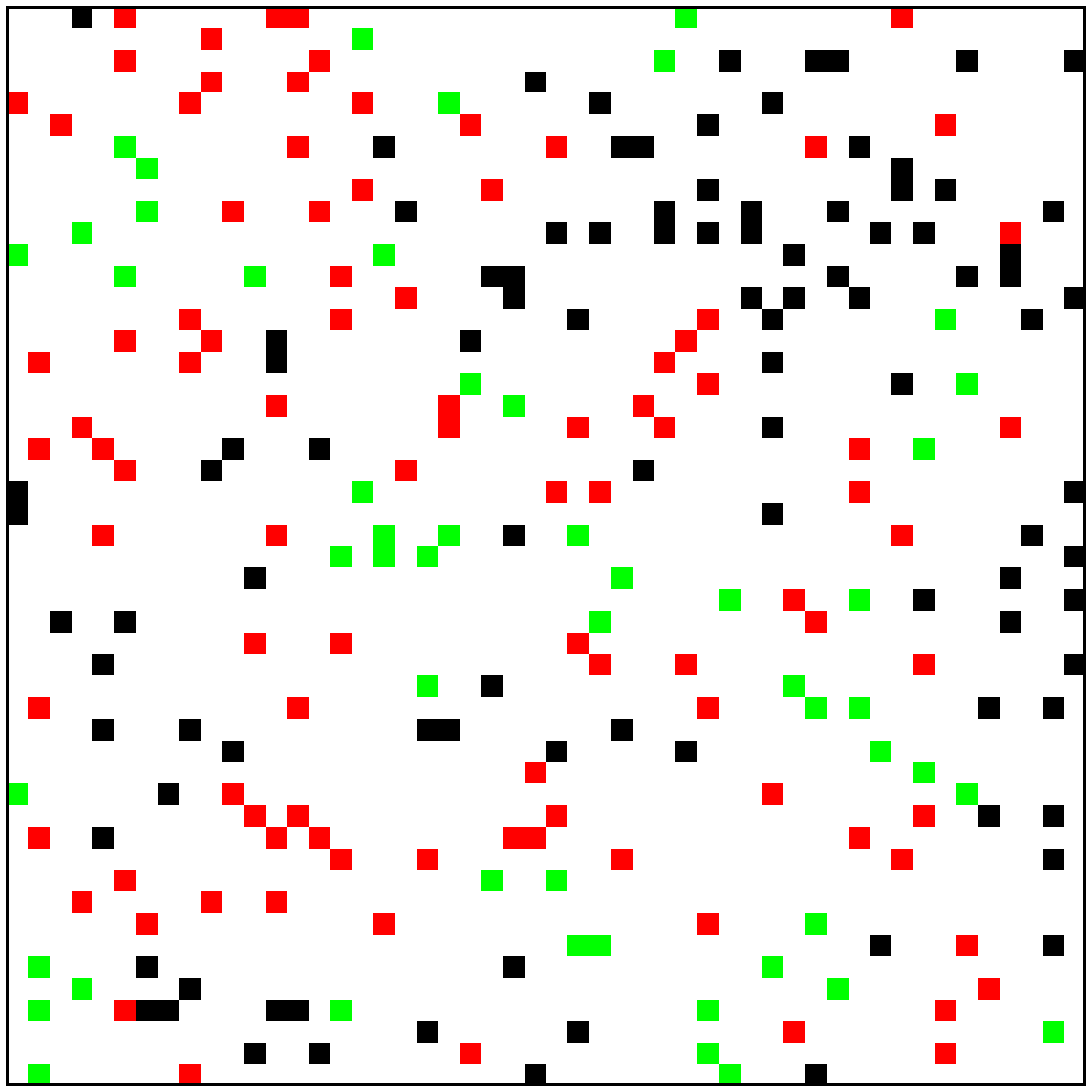}
\includegraphics[scale=0.22,angle=-90,draft=false]{\figs 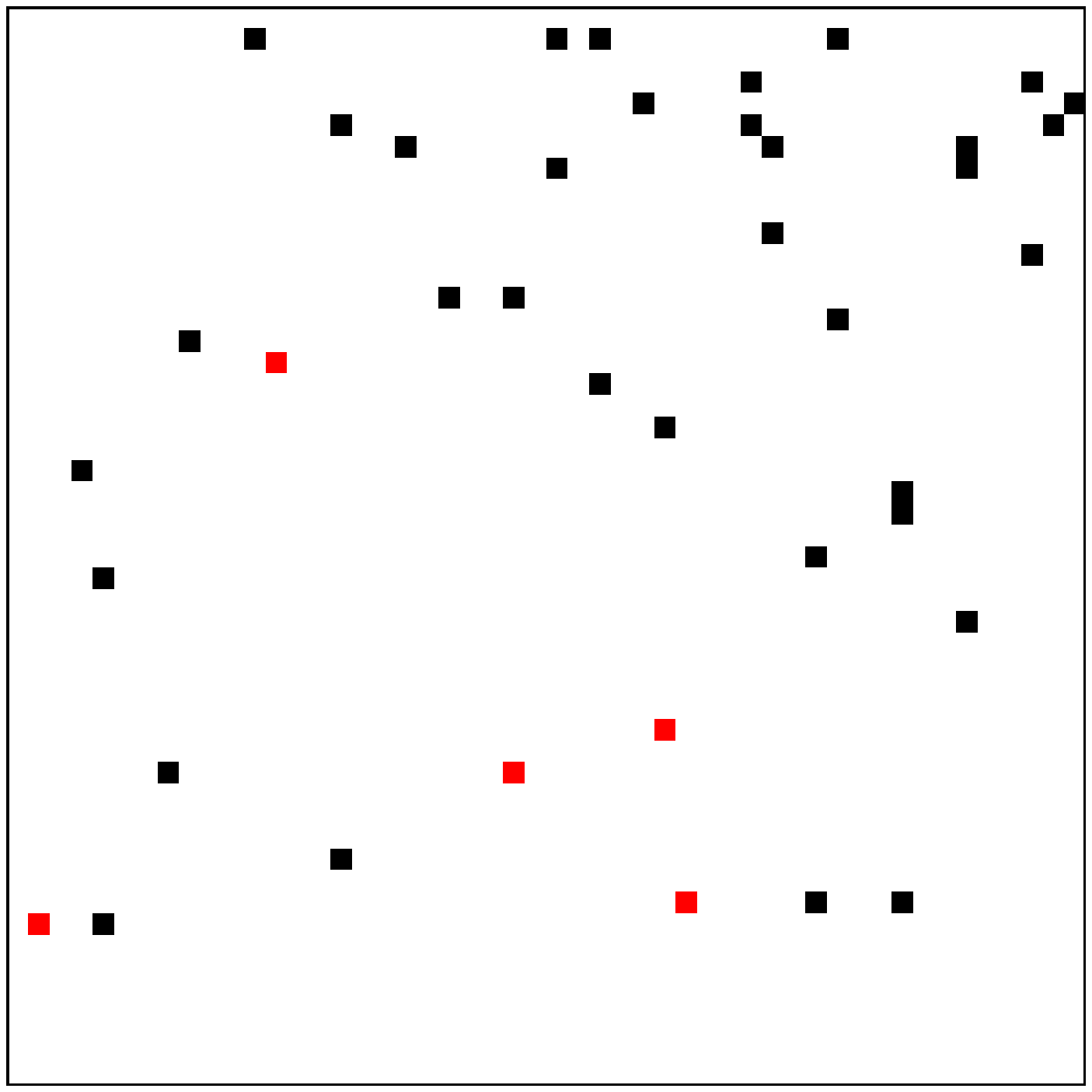}
\end{center}
\caption{Snapshots corresponding to different stages of a simulation carried out for the random truel with initial proportions of  $x_A=0.3$ (black), $x_B=0.3$ (red) and $x_C=0.4$ (green), for a set of $N=2500$ players arranged in a two--dimensional grid. The marksmanships are $a=1$, $b=0.8$, $c=0.5$.\label{evolution_random}}
\end{figure}

In this collective truel,  the group that will survive at the end depends, for a fixed values of the marksmanships, on the initial proportions of players. This dependence is summarized in Fig.~\ref{diagram_random}, where we plot in a color code the group that has the highest winning probability as a functions of the  initial proportions.

\begin{figure}[hbt!]
\begin{center}
\includegraphics[scale=0.4,angle=-90,draft=false]{\figs 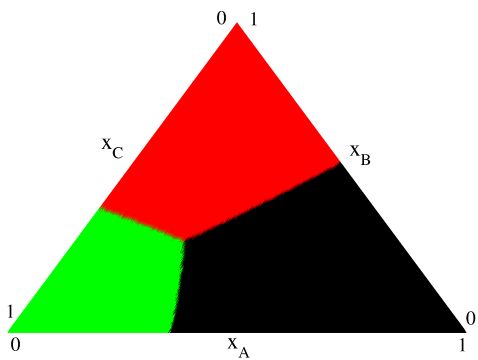}
\end{center}
\caption{Diagram of  the favorite group, the one with the highest probability of winning, (black corresponds to group A, red to  B and green to C) in terms of the initial proportions $x_A$, $x_B$ and $x_C$, for a set of $N=400$ players arranged in a two--dimensional grid and playing the random collective  truel. The probabilities have been obtained after a large number of simulations.\label{diagram_random}}
\end{figure}

It is easy to modify step 2 by considering the rule of the sequential truel by which players shoot in inverse order to their marksmanship. 
A typical realization is shown in Fig.\ref{evolution_sequential}. In this occasion the winning group is the weakest one, group C. This \textit{survival of the weakest} effect is also present in the diagram of Fig.~\ref{diagram_sequential}, as now groups B and C have increased the region in parameter space where they win the truel, compared to the diagram of the random truel in Fig.\ref{diagram_random}.

\begin{figure}[hbt!]
\begin{center}
\includegraphics[scale=0.22,angle=-90,draft=false]{\figs 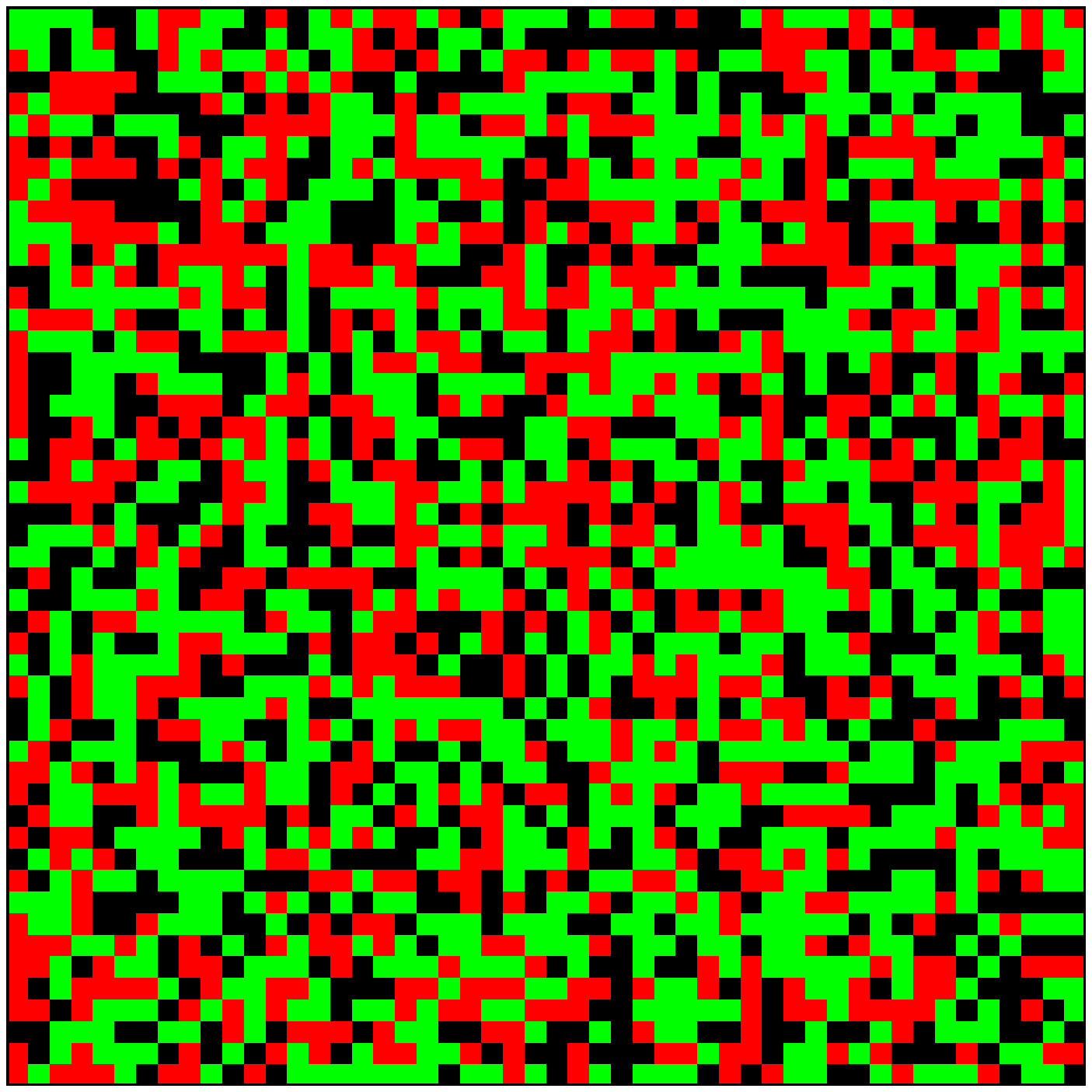}
\includegraphics[scale=0.22,angle=-90,draft=false]{\figs 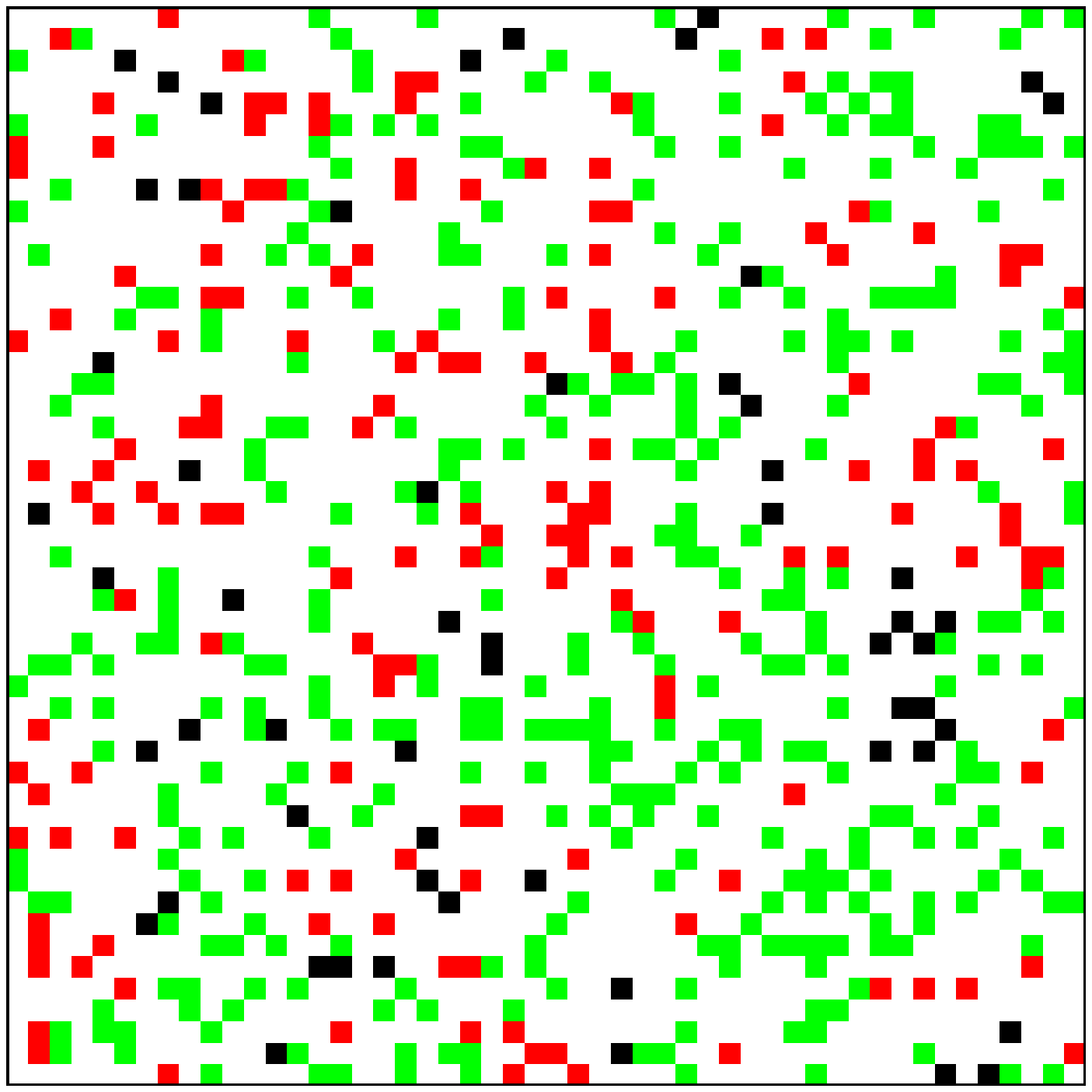}
\includegraphics[scale=0.22,angle=-90,draft=false]{\figs 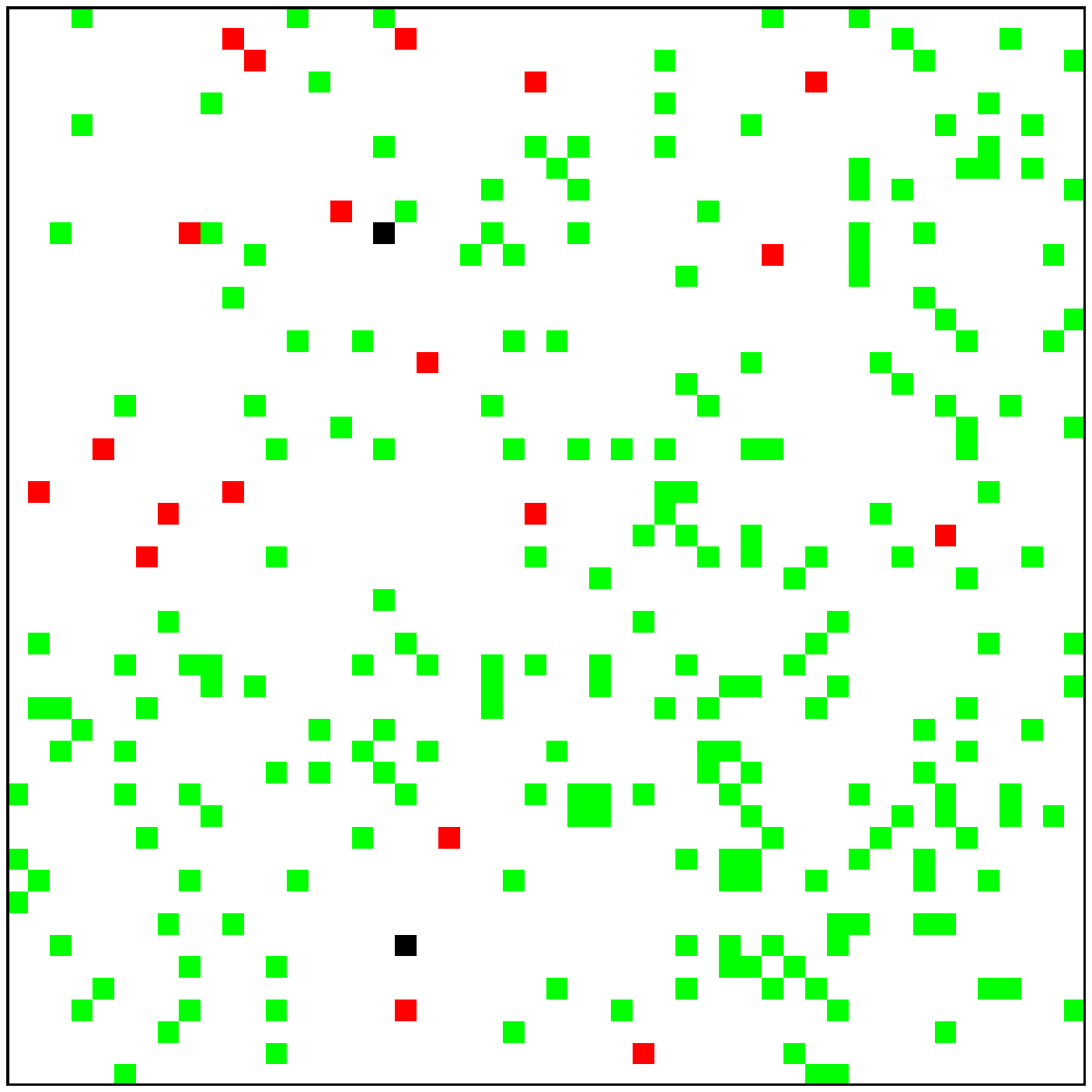}
\includegraphics[scale=0.22,angle=-90,draft=false]{\figs 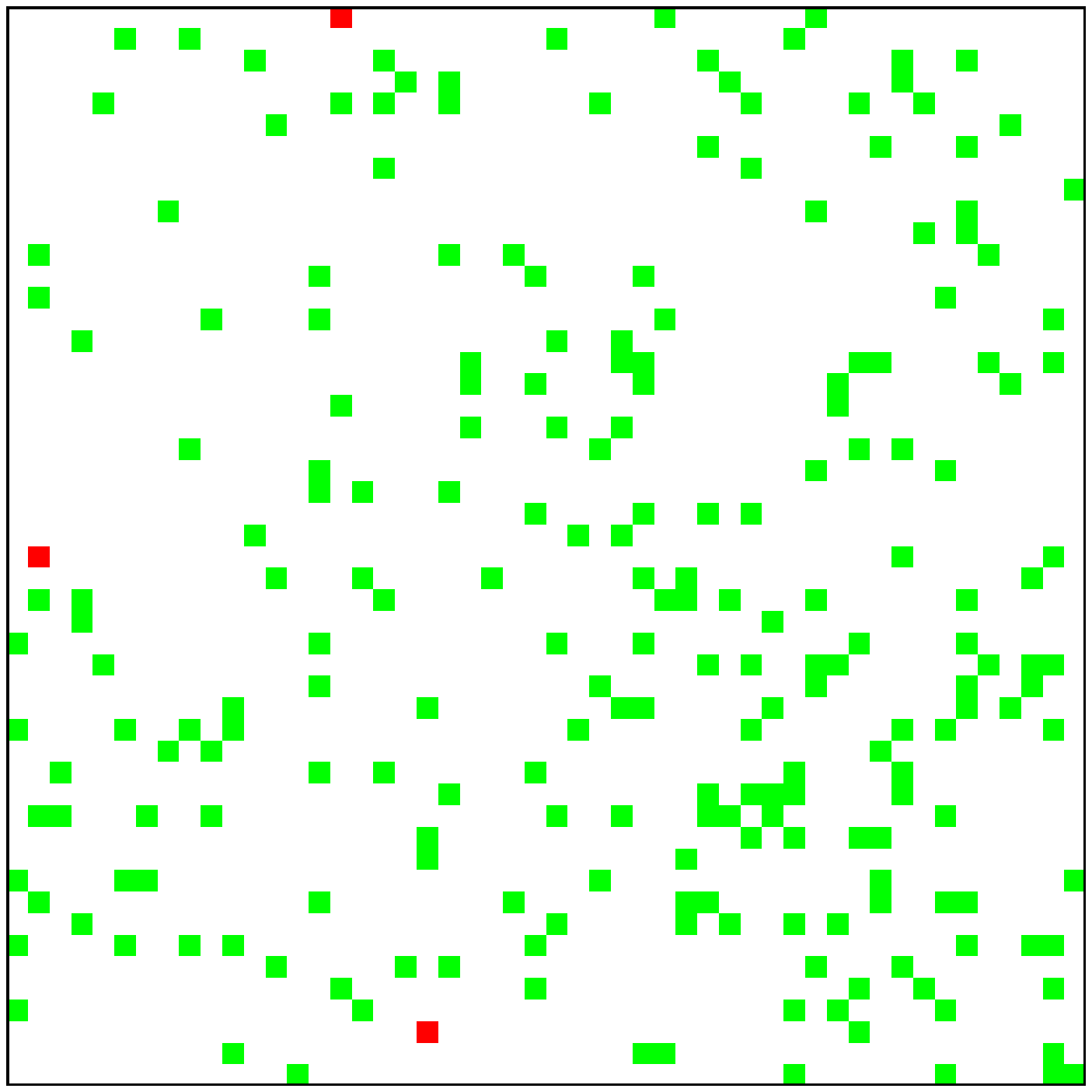}
\end{center}
\caption{Same as Fig.\ref{evolution_random} but using the rules of the sequential truel.\label{evolution_sequential}}
\end{figure}

It is possible to distinguish two different regimes in the dynamics. Almost all truel competitions take place during the first steps where a large fraction of the population is removed. At the end of this first regime, the largest remaining population is the one that possesses the higher survival probability when playing a single truel and the system presents many empty sites. Later, in a second regime, players start to diffuse to neighboring sites increasing the appearance of duel encounters. Consequently, the evolution will result from a balance between the population favored by the existence of \textit{duels} (the one with the highest marksmanship), and the one favored by possessing a high proportion of the remaining population.

\begin{figure}[hbt!]
\begin{center}
\includegraphics[scale=0.4,angle=-90,draft=false]{\figs 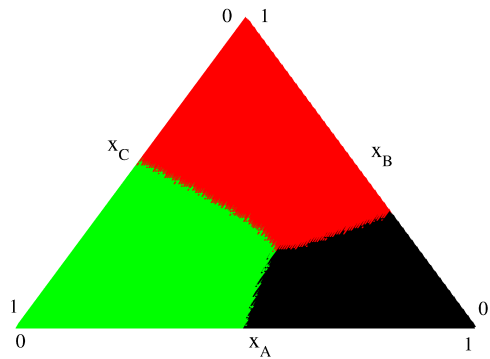}
\end{center}
\caption{Same as Fig.\ref{diagram_random} but using the rules of the sequential truel.\label{diagram_sequential}}
\end{figure}

Finally, in Fig.~\ref{evolution_opinion} we show some snapshots of a simulation carried out for the case of the opinion model. For the set of marksmanships chosen $a=1$, $b=0.8$ and $c=0.5$ we find the favorite opinion is always the one with highest marksmanship, A. This occurs even for very small initial proportion $x_A$ and it is a reflection of the large region in parameter space where A becomes the favorite opinion, as it was shown in Fig.\ref{opinion}.

\begin{figure}[hbt!]
\begin{center}
\includegraphics[scale=0.22,angle=-90,draft=false]{\figs 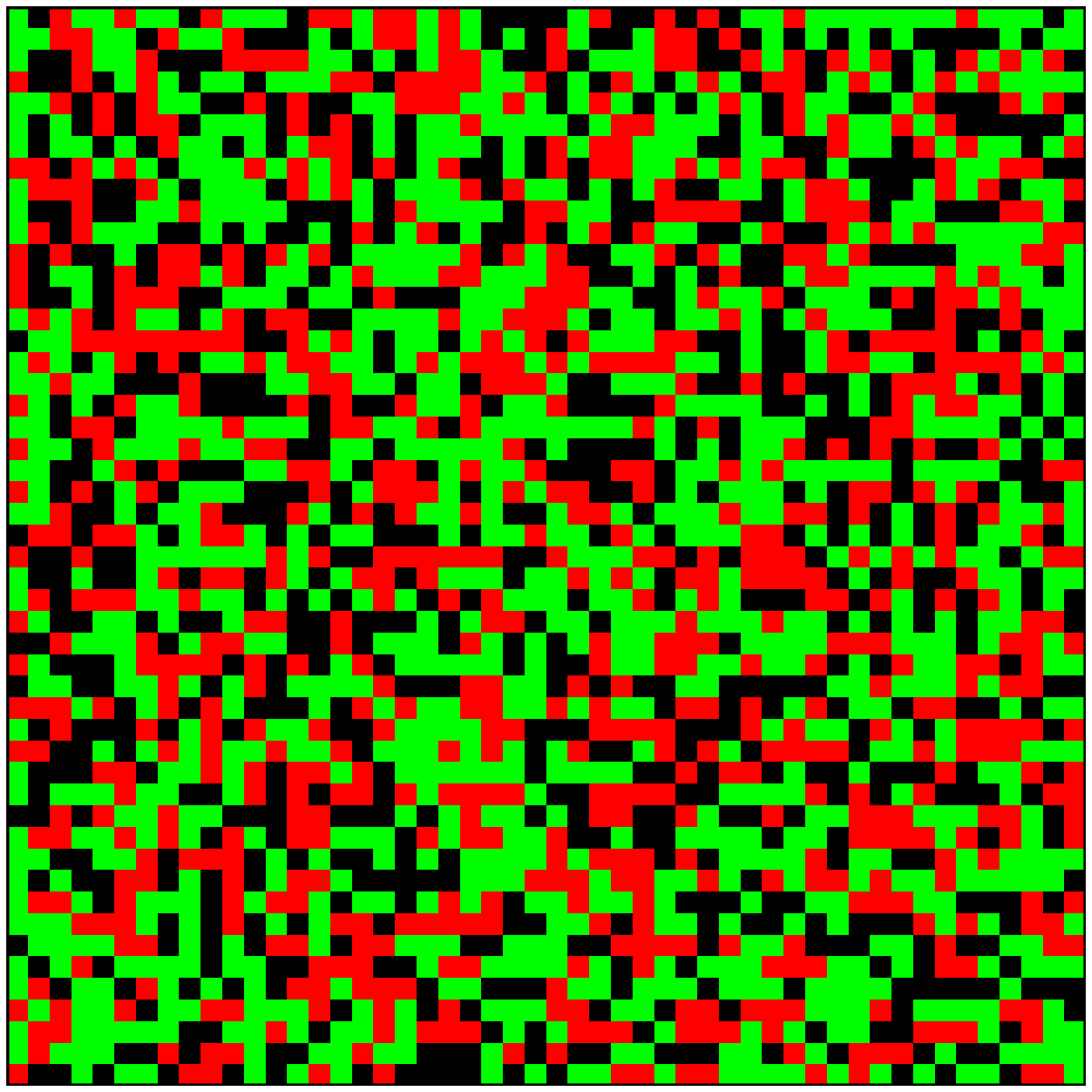}
\includegraphics[scale=0.22,angle=-90,draft=false]{\figs 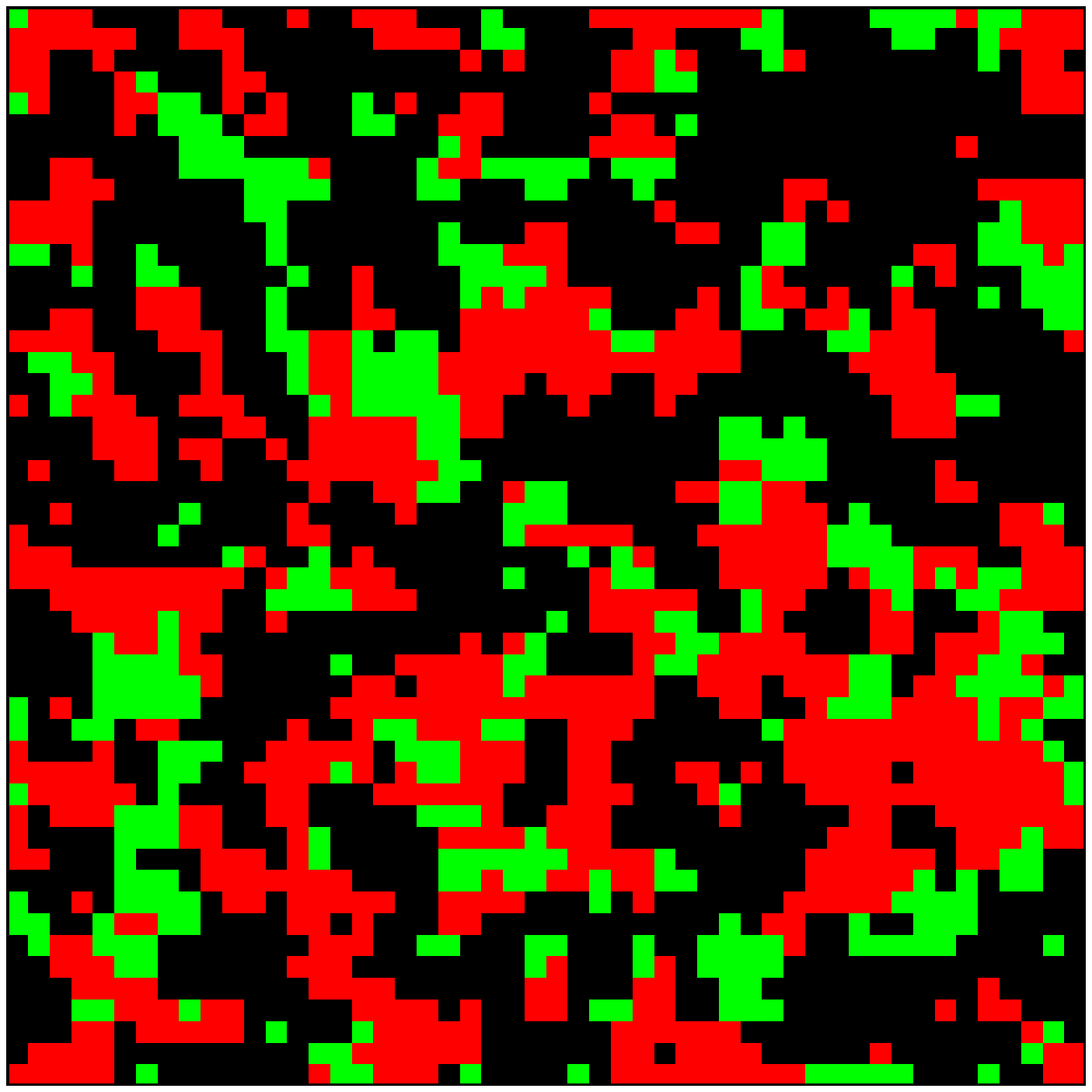}
\includegraphics[scale=0.22,angle=-90,draft=false]{\figs 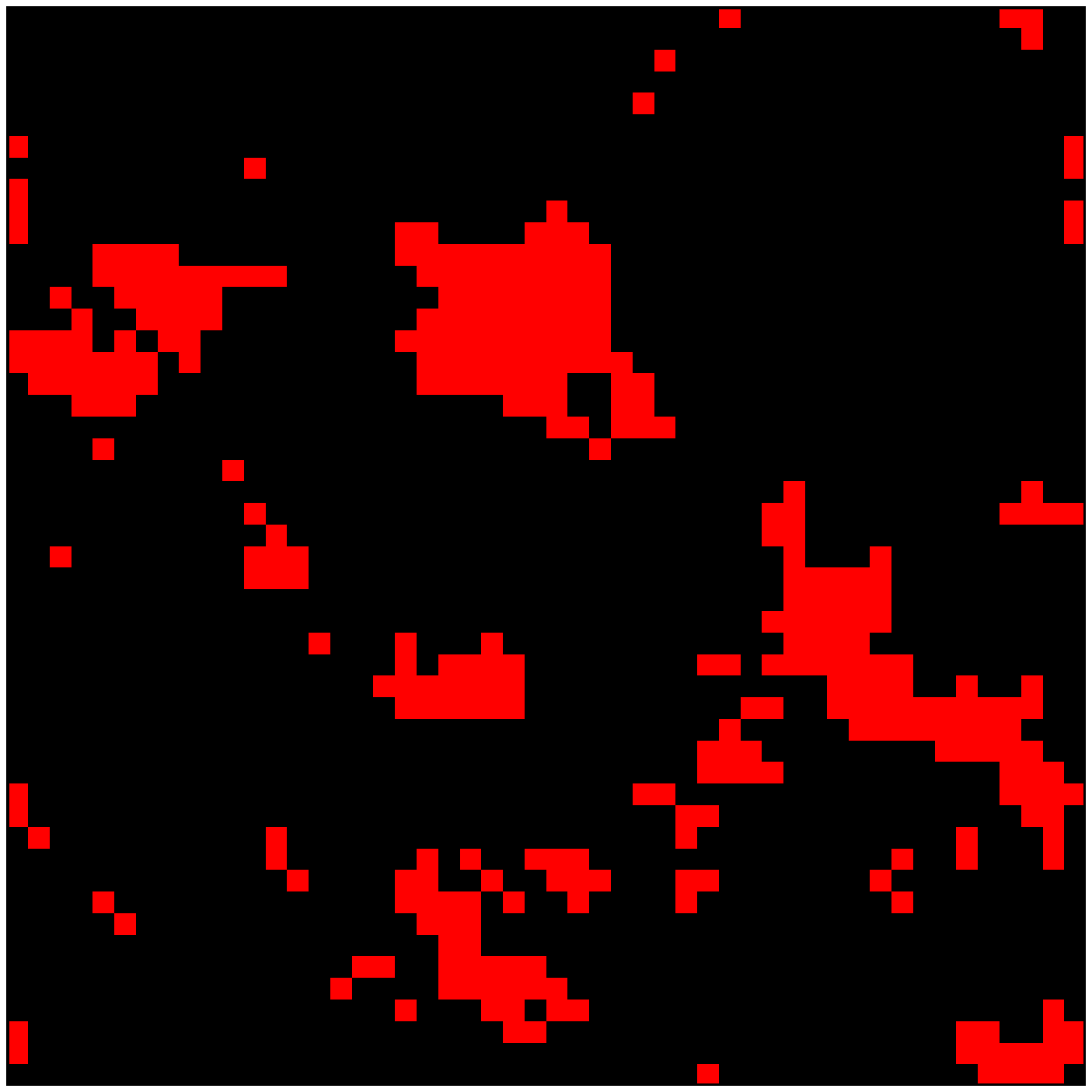}
\includegraphics[scale=0.22,angle=-90,draft=false]{\figs 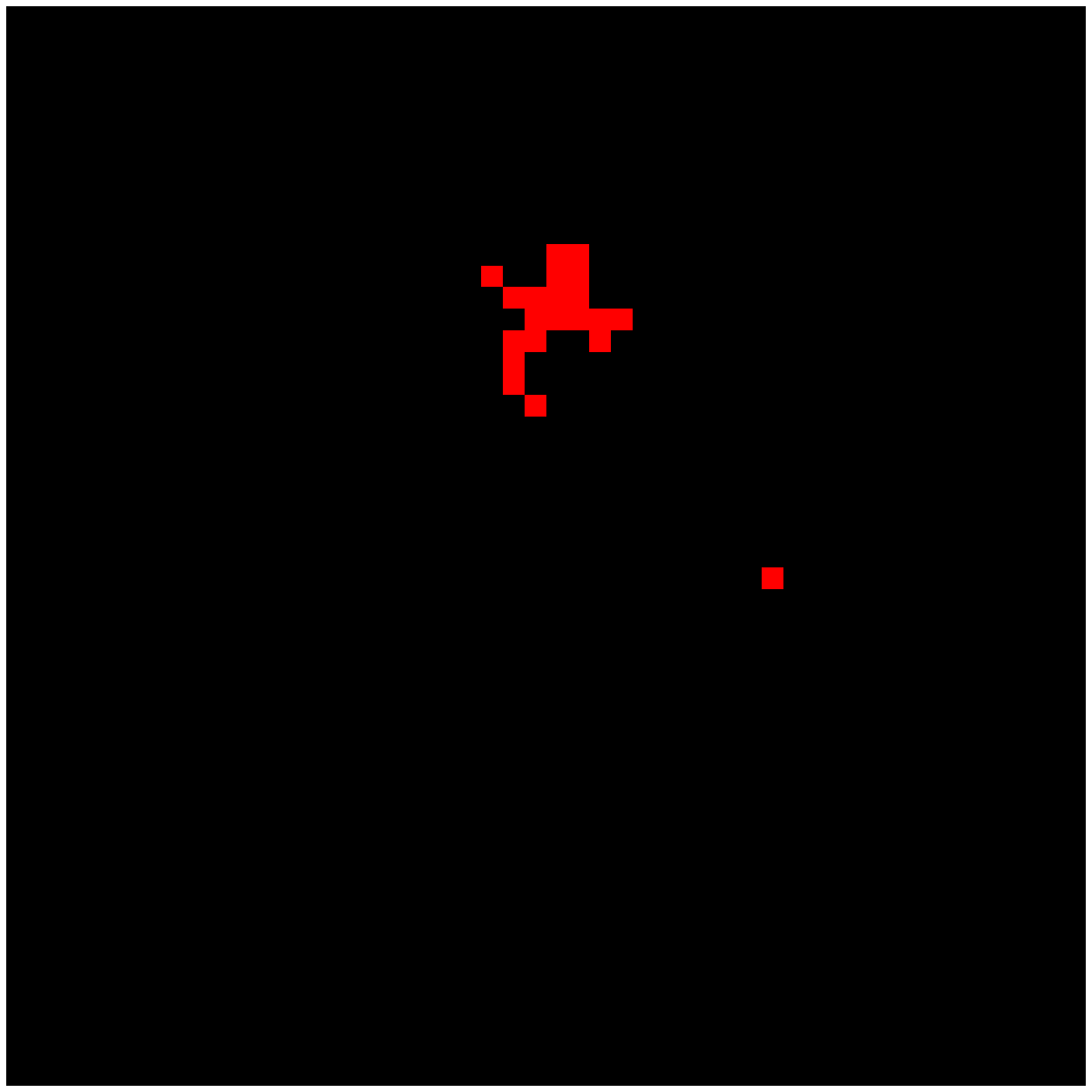}
\end{center}
\caption{Same as Fig.\ref{evolution_random} but using the rules of the opinion truel. \label{evolution_opinion}}
\end{figure}

\section{Summary}

We have reviewed in this work some counterintuitive results concerning truel games, an extension of a duel between three people. The truels were originally analyzed using the tools of game theory and provide a clarifying  example of the concept of Nash equilibrium. The paradox is that the use of the best possible strategy implies that the {\sl a priori}�worst performing player might be the favorite to win the game, something that we could name as the {\sl failure of the fittest}. This result conforms to the observation that in some situations the most appropriate candidate does not succeed because of competition with other less suited individuals.  We have recently revisited those games from the point of view of stochastic processes, a technique that allows us to reproduce some of the known results for the equilibrium points as, for instance, the regions in which every kind of player is favorite to win the game. Besides showing results for the traditional random and sequential versions of the truel, we have introduced a simple modification that allows us to interpret this game in terms of a model for opinion spreading.

We have addressed the question of the distribution of winners in a truel league. We have shown that under certain circumstances is not always recommendable to be the best player. Indeed, when playing the sequential truel intermediate values of the marksmanship.
perform better on the average.  We have also generalized the random truel for a number of players greater than three. In this case, already for $4$ players it is noticeable the appearance of an optimum value for the marksmanship lower than one. As the number of players increases, the optimum value shifts towards lower and lower values of the marksmanship. This indicates that the paradoxical effect mentioned above increases with the number of players.

Finally, we have introduced a spatial structure in the set of players and their connectivities by simulating the different versions of the truels in a two--dimensional grid.   We have shown the importance of the dynamics in the early stages of the simulation. In this case, the winner  of the game depends on the initial population of players belonging to each possible class.

\bigskip

We acknowledge financial support from the Spanish Government and FEDER (EU) through projects FIS2004-5073, FIS2004-953. P.A. is supported by a grant from the government of the Balearic Islands.

\bibliographystyle{unsrt}
\bibliography{\dir Bibliography}

\end{document}